\documentclass{amsart}
\usepackage{amssymb}

\usepackage{amscd}
\usepackage{thmdefs}

\input{tcilatex}
\theoremstyle{definition}
\theoremstyle{remark}
\numberwithin{equation}{section}

\input tcilatex

\begin{document}
\title[Amalgamated R-Diagonal Pairs]{Amalgamated R-Diagonal Pairs}
\author{Ilwoo Cho}
\address{Dep. of Math, Univ. of Iowa, Iowa City, IA, U.S.A}
\email{ilcho@math.uiowa.edu}
\maketitle

\begin{abstract}
In this paper, we will consider the properties of amalgamated R-diagonal
pairs. We characterize the amalgamated R-diagonality of pairs of amalgamated
random variables by certain cumulant-relation.
\end{abstract}

\strut

Voiculescu developed Free Probability Theory. Here, the classical concept of
Independence in Probability theory is replaced by a noncommutative analogue
called Freeness (See [7]). There are two approaches to study Free
Probability Theory. One of them is the original analytic approach of
Voiculescu (See [7] and [10]) and the other one is the combinatorial
approach of Speicher and Nica (See [10], [1] and [11]). Speicher defined the
free cumulants which are the main objects in the combinatorial approach of
Free Probability Theory. And he developed free probability theory by using
the combinatorics and lattice theory on collections of noncrossing
partitions (See [11]). Also, Speicher considered the operator-valued free
probability theory, which is also defined and observed analytically by
Voiculescu, when $\Bbb{C}$ is replaced to an arbitrary algebra $B$ (See [10]
and [10]). Nica defined R-transforms of several random variables (See [1]).
He defined these R-transforms as multivariable formal series in
noncommutative several indeterminants. To observe the R-transform, the M\"{o}%
bius Inversion under the embedding of lattices plays a key role (See
[10],[11],[3],[8],[9] and [17]). In this paper, we will consider the $B$%
-even elements and R-diagonal pairs of $B$-valued random variables. Let $%
(A,\varphi )$ be a NCPSpace over $B$ and let $x_{1},x_{2}\in (A,\varphi )$
be $B$-valued random variables. We say that a pair $(x_{1},x_{2})$ is
R-diagonal if there exists a $B$-formal series $f_{1},\,g\in \Theta _{B}^{1}$
such that

\strut 

\begin{center}
$R_{x_{1},x_{2}}(z_{1},z_{2})=f(z_{1}z_{2})+g(z_{2}z_{1}).$
\end{center}

\strut 

We call this $B$-formal series $\ (f,$ $g)\in \Theta _{B}^{1}\times \Theta
_{B}^{1}$ \ the determining series of the pair $(x_{1},x_{2}).$ We show that
the determining series $\left( f,g\right) $ is determined by

\strut 

\begin{center}
$f=R_{xy}\,\,\frame{*}_{B}\,\,Mob$ \ \ and \ \ $g=R_{yx}\,\,\frame{*}%
_{B}\,\,\,Mob.$
\end{center}

\strut \strut 

Also, similar to the scalar-valued case observed by Nica and Speicher, we
can get the following fact that if the random variables $a,a^{\prime }\in
(A,\varphi )$ $B$-even and if they are free over $B,$ then the pair $%
(aa^{\prime },a^{\prime }a)$ is $B$-valued R-diagonal.

\strut

\strut

\section{Amalgamated R-transform Theory}

\strut

\strut

In this section, we will define an R-transform of several $B$-valued random
variables. Note that to study R-transforms is to study operator-valued
distributions. R-transforms with single variable is defined by Voiculescu
(over $B,$ in particular, $B=\Bbb{C}$. See [7] and [10]). Over $\Bbb{C},$
Nica defined multi-variable R-transforms in [1]. In [8], we extended his
concepts, over $B.$ R-transforms of $B$-valued random variables can be
defined as $B$-formal series with its $(i_{1},...,i_{n})$-th coefficients, $%
(i_{1},...,i_{n})$-th cumulants of $B$-valued random variables, where $%
(i_{1},...,i_{n})\in \{1,...,s\}^{n},$ $\forall n\in \Bbb{N}.$

\strut

\begin{definition}
Let $(A,\varphi )$ be a NCPSpace over $B$ and let $x_{1},...,x_{s}\in
(A,\varphi )$ be $B$-valued random variables ($s\in \Bbb{N}$). Let $%
z_{1},...,z_{s}$ be noncommutative indeterminants. Define a moment series of 
$x_{1},...,x_{s}$, as a $B$-formal series, by

\strut 

\begin{center}
$M_{x_{1},...,x_{s}}(z_{1},...,z_{s})=\sum_{n=1}^{\infty }\underset{%
i_{1},..,i_{n}\in \{1,...,s\}}{\sum }\varphi
(x_{i_{1}}b_{i_{2}}x_{i_{2}}...b_{i_{n}}x_{i_{n}})\,z_{i_{1}}...z_{i_{n}},$
\end{center}

\strut 

where $b_{i_{2}},...,b_{i_{n}}\in B$ are arbitrary for all $%
(i_{2},...,i_{n})\in \{1,...,s\}^{n-1},$ $\forall n\in \Bbb{N}.$

\strut 

Define an R-transform of $x_{1},...,x_{s}$, as a $B$-formal series, by

\strut 

\begin{center}
$R_{x_{1},...,x_{s}}(z_{1},...,z_{s})=\sum_{n=1}^{\infty }\underset{%
i_{1},...,i_{n}\in \{1,...,s\}}{\sum }k_{n}(x_{i_{1}},...,x_{i_{n}})%
\,z_{i_{1}}...z_{i_{n}},$
\end{center}

\strut with

\begin{center}
$k_{n}(x_{i_{1}},...,x_{i_{n}})=c^{(n)}(x_{i_{1}}\otimes
b_{i_{2}}x_{i_{2}}\otimes ...\otimes b_{i_{n}}x_{i_{n}}),$
\end{center}

\strut 

where $b_{i_{2}},...,b_{i_{n}}\in B$ are arbitrary for all $%
(i_{2},...,i_{n})\in \{1,...,s\}^{n-1},$ $\forall n\in \Bbb{N}.$ Here, $%
\widehat{c}=(c^{(n)})_{n=1}^{\infty }$ is a cumulant multiplicative function
induced by $\varphi $ in $I(A,B).$
\end{definition}

\strut

Denote a set of all $B$-formal series with $s$-noncommutative indeterminants
($s\in \Bbb{N}$), by $\Theta _{B}^{s}$. i.e if $g\in \Theta _{B}^{s},$ then

\begin{center}
$g(z_{1},...,z_{s})=\sum_{n=1}^{\infty }\underset{i_{1},...,i_{n}\in
\{1,...,s\}}{\sum }b_{i_{1},...,i_{n}}\,z_{i_{1}}...z_{i_{n}},$
\end{center}

\strut

where $b_{i_{1},...,i_{n}}\in B,$ for all $(i_{1},...,i_{n})\in
\{1,...,s\}^{n},$ $\forall n\in \Bbb{N}.$ Trivially, by definition, $%
M_{x_{1},...,x_{s}},$ $R_{x_{1},...,x_{s}}\in \Theta _{B}^{s}.$ By $\mathcal{%
R}_{B}^{s},$\ we denote a set of all R-transforms of $s$-$B$-valued random
variables. Recall that, set-theoratically,

\begin{center}
$\Theta _{B}^{s}=\mathcal{R}_{B}^{s},$ sor all $s\in \Bbb{N}.$
\end{center}

\strut

We can also define symmetric moment series and symmetric R-transform by $%
b_{0}\in B,$ by

\strut

\begin{center}
$M_{x_{1},...,x_{s}}^{symm(b_{0})}(z_{1},...,z_{s})=\sum_{n=1}^{\infty }%
\underset{i_{1},...,i_{n}\in \{1,...,s\}}{\sum }\varphi
(x_{i_{1}}b_{0}x_{i_{2}}...b_{0}x_{i_{n}})\,z_{i_{1}}...z_{i_{n}}$
\end{center}

and

\begin{center}
$R_{x_{1},...,x_{s}}^{symm(b_{0})}(z_{1},...,z_{s})=\sum_{n=1}^{\infty }%
\underset{i_{1},..,i_{n}\in \{1,...,s\}}{\sum }%
k_{n}^{symm(b_{0})}(x_{i_{1}},...,x_{i_{n}})\,z_{i_{1}}...z_{i_{n}},$
\end{center}

with

\begin{center}
$k_{n}^{symm(b_{0})}(x_{i_{1}},...,x_{i_{n}})=c^{(n)}(x_{i_{1}}\otimes
b_{0}x_{i_{2}}\otimes ...\otimes b_{0}x_{i_{n}}),$
\end{center}

\strut \strut

for all $(i_{1},...,i_{n})\in \{1,...,s\}^{n},$ $\forall n\in \Bbb{N}.$

\strut

If $b_{0}=1_{B},$ then we have trivial moment series and trivial R-transform
of $x_{1},...,x_{s}$ denoted by $M_{x_{1},...,x_{s}}^{t}$ and $%
R_{x_{1},...,x_{s}}^{t},$ respectively. By definition, for the fixed random
variables $x_{1},...,x_{s}\in (A,\varphi ),$ there are infinitely many
R-transforms of them (resp. moment series of them). Symmetric and trivial
R-transforms of them are special examples. Let

\strut 

\begin{center}
$C=\underset{(i_{1},...,i_{n})\in \Bbb{N}^{n}}{\cup }%
\{(1_{B},b_{i_{2}},...,b_{i_{n}}):b_{i_{j}}\in B\}.$
\end{center}

\strut 

Suppose that we have

\strut 

\begin{center}
$coef_{i_{1},...,i_{n}}\left( R_{x_{1},...,x_{s}}\right) =c^{(n)}\left(
x_{i_{1}}\otimes b_{i_{2}}x_{i_{2}}\otimes ...\otimes
b_{i_{n}}x_{i_{n}}\right) ,$
\end{center}

\strut 

where $(1_{B},b_{i_{2}},...,b_{i_{n}})\in C,$ for all $(i_{1},...,i_{n})\in 
\Bbb{N}^{n}.$ Then we can rewite the R-transform of $x_{1},...,x_{s},$ $%
R_{x_{1},...,x_{s}}$ by $R_{x_{1},...,x_{s}}^{C}.$ If $C_{1}$ and $C_{2}$
are such collections, then in general $R_{x_{1},...,x_{s}}^{C_{1}}\neq
R_{x_{1},...,x_{s}}^{C_{2}}$ (resp. $M_{x_{1},...,x_{s}}^{C_{1}}\neq
M_{x_{1},...,x_{s}}^{C_{2}}$). From now, for the random variables $%
x_{1},...,x_{s},y_{1},...,y_{s},$ if we write $R_{x_{1},...,x_{s}}$ and $%
R_{y_{1},...,y_{s}},$ then it means that $%
R_{x_{1},...,x_{s}}^{C}=R_{y_{1},....,y_{s}}^{C},$ for the same collection $%
C.$ If there's no confusion, we will omit to write such collection. The
followings are known in [10] and [8] ;

\strut

\begin{proposition}
Let $(A,\varphi )$ be a NCPSpace over $B$ and let $%
x_{1},...,x_{s},y_{1},...,y_{p}\in (A,\varphi )$ be $B$-valued random
variables, where $s,p\in \Bbb{N}.$ Suppose that $\{x_{1},...,x_{s}\}$ and $%
\{y_{1},...,y_{p}\}$ are free in $(A,\varphi ).$ Then

\strut 

(1) $
R_{x_{1},...,x_{s},y_{1},...,y_{p}}(z_{1},...,z_{s+p})=R_{x_{1},...,x_{s}}(z_{1},...,z_{s})+R_{y_{1},...,y_{p}}(z_{1},...,z_{p}).
$

\strut 

(2) If $s=p,$ then $R_{x_{1}+y_{1},...,x_{s}+y_{s}}(z_{1},...,z_{s})=\left(
R_{x_{1},...,x_{s}}+R_{y_{1},...,y_{s}}\right) (z_{1},...,z_{s}).$

$\square $
\end{proposition}

\strut \strut

Note that if $f,g\in \Theta _{B}^{s},$ then we can always choose free $%
\{x_{1},...,x_{s}\}$ and $\{y_{1},...,y_{s}\}$ in (some) NCPSpace over $B,$ $%
(A,\varphi ),$ such that

\begin{center}
$f=R_{x_{1},...,x_{s}}$ \ \ and \ \ $g=R_{y_{1},...,y_{s}}.$
\end{center}

\strut

\begin{definition}
(1) Let $s\in \Bbb{N}.$ Let $(f,g)\in \Theta _{B}^{s}\times \Theta _{B}^{s}.$
Define \frame{*}\thinspace \thinspace $:\Theta _{B}^{s}\times \Theta
_{B}^{s}\rightarrow \Theta _{B}^{s}$ by

\strut 

\begin{center}
$\left( f,g\right) =\left(
R_{x_{1},...,x_{s}}^{C_{1}},\,R_{y_{1},...,y_{s}}^{C_{2}}\right) \longmapsto
R_{x_{1},...,x_{s}}^{C_{1}}\,\,\frame{*}\,\,R_{y_{1},...,y_{s}}^{C_{2}}.$
\end{center}

\strut 

Here, $\{x_{1},...,x_{s}\}$ and $\{y_{1},...,y_{s}\}$ are free in $%
(A,\varphi )$. Suppose that

\strut 

\begin{center}
$coef_{i_{1},..,i_{n}}\left( R_{x_{1},...,x_{s}}^{C_{1}}\right)
=c^{(n)}(x_{i_{1}}\otimes b_{i_{2}}x_{i_{2}}\otimes ...\otimes
b_{i_{n}}x_{i_{n}})$
\end{center}

and

\begin{center}
$coef_{i_{1},...,i_{n}}(R_{y_{1},...,y_{s}}^{C_{2}})=c^{(n)}(y_{i_{1}}%
\otimes b_{i_{2}}^{\prime }y_{i_{2}}\otimes ...\otimes b_{i_{n}}^{\prime
}y_{i_{n}}),$
\end{center}

\strut 

for all $(i_{1},...,i_{n})\in \{1,...,s\}^{n},$ $n\in \Bbb{N},$ where $%
b_{i_{j}},b_{i_{n}}^{\prime }\in B$ arbitrary. Then

\strut 

$coef_{i_{1},...,i_{n}}\left( R_{x_{1},...,x_{s}}^{C_{1}}\,\,\frame{*}%
\,\,R_{y_{1},...,y_{s}}^{C_{2}}\right) $

\strut 

$=\underset{\pi \in NC(n)}{\sum }\left( \widehat{c_{x}}\oplus \widehat{c_{y}}%
\right) (\pi \cup Kr(\pi ))(x_{i_{1}}\otimes y_{i_{1}}\otimes
b_{i_{2}}x_{i_{2}}\otimes b_{i_{2}}^{\prime }y_{i_{2}}\otimes ...\otimes
b_{i_{n}}x_{i_{n}}\otimes b_{i_{n}}^{\prime }y_{i_{n}})$

$\strut $

$\overset{denote}{=}\underset{\pi \in NC(n)}{\sum }\left( k_{\pi
}^{C_{1}}\oplus k_{Kr(\pi )}^{C_{2}}\right)
(x_{i_{1}},y_{i_{1}},...,x_{i_{n}}y_{i_{n}}),$

\strut \strut 

where $\widehat{c_{x}}\oplus \widehat{c_{y}}=\widehat{c}\mid
_{A_{x}*_{B}A_{y}},$ $A_{x}=A\lg \left( \{x_{i}\}_{i=1}^{s},B\right) $ and $%
A_{y}=A\lg \left( \{y_{i}\}_{i=1}^{s},B\right) $ and where $\pi \cup Kr(\pi )
$ is an alternating union of partitions in $NC(2n)$
\end{definition}

\strut

\begin{proposition}
(See [8])\strut Let $(A,\varphi )$ be a NCPSpace over $B$ and let $%
x_{1},...,x_{s},y_{1},...,y_{s}\in (A,\varphi )$ be $B$-valued random
variables ($s\in \Bbb{N}$). If $\{x_{1},...,x_{s}\}$ and $\{y_{1},...,y_{s}\}
$ are free in $(A,\varphi ),$ then we have

\strut 

$k_{n}(x_{i_{1}}y_{i_{1}},...,x_{i_{n}}y_{i_{n}})$

\strut 

$=\underset{\pi \in NC(n)}{\sum }\left( \widehat{c_{x}}\oplus \widehat{c_{y}}%
\right) (\pi \cup Kr(\pi ))(x_{i_{1}}\otimes y_{i_{1}}\otimes
b_{i_{2}}x_{i_{2}}\otimes y_{i_{2}}\otimes ...\otimes
b_{i_{n}}x_{i_{n}}\otimes y_{i_{n}})$

\strut 

$\overset{denote}{=}\underset{\pi \in NC(n)}{\sum }\left( k_{\pi }\oplus
k_{Kr(\pi )}^{symm(1_{B})}\right)
(x_{i_{1}},y_{i_{1}},...,x_{i_{n}},y_{i_{n}}),$

\strut 

for all $(i_{1},...,i_{n})\in \{1,...,s\}^{n},$ $\forall n\in \Bbb{N},$ $%
b_{i_{2}},...,b_{i_{n}}\in B,$ arbitrary, where $\widehat{c_{x}}\oplus 
\widehat{c_{y}}=\widehat{c}\mid _{A_{x}*_{B}A_{y}},$ $A_{x}=A\lg \left(
\{x_{i}\}_{i=1}^{s},B\right) $ and $A_{y}=A\lg \left(
\{y_{i}\}_{i=1}^{s},B\right) .$ \ $\square $
\end{proposition}

\strut

This shows that ;

\strut

\begin{corollary}
(See [8]) Under the same condition with the previous proposition,

\strut 

\begin{center}
$R_{x_{1},...,x_{s}}\,\,\frame{*}\,%
\,R_{y_{1},...,y_{s}}^{t}=R_{x_{1}y_{1},...,x_{s}y_{s}}.$
\end{center}

$\square $
\end{corollary}

\strut

Notice that, in general, unless $b_{i_{2}}^{\prime }=...=b_{i_{n}}^{\prime
}=1_{B}$ in $B,$

\strut

\begin{center}
$R_{x_{1},...,x_{s}}^{C_{1}}\,\,\frame{*}\,\,R_{y_{1},...,y_{s}}^{C_{2}}\neq
R_{x_{1}y_{1},...,x_{s}y_{s}}^{C_{j}},$ \ \ $j=1.2.$
\end{center}

\strut

However, as we can see above,

\strut

\begin{center}
$R_{x_{1},...,x_{s}}\,\,\frame{*}\,%
\,R_{y_{1},...,y_{s}}^{t}=R_{x_{1}y_{1},...,x_{s}y_{s}}$
\end{center}

and

\begin{center}
$R_{x_{1},...,x_{s}}^{t}\,\,\frame{*}\,%
\,R_{y_{1},...,y_{s}}^{t}=R_{x_{1}y_{1},...,x_{s}y_{s}}^{t},$
\end{center}

\strut

where $\{x_{1},...,x_{s}\}$ and $\{y_{1},...,y_{s}\}$ are free over $B.$
Over $B=\Bbb{C},$ the last equation is proved by Nica and Speicher in [1]
and [11]. Actually, their R-transforms (over $\Bbb{C}$) is our trivial
R-transforms (over $\Bbb{C}$).

\strut

\strut

\strut

\section{R-diagonal Pairs}

\strut

\strut

\subsection{$B$-valued Even Random Variables}

\strut

\strut

\strut In this section, we will consider the $B$-evenness. Let $(A,\varphi )$
be a NCPSpace over $B,$ with its $B$-trace $\varphi :A\rightarrow B.$

\strut

\begin{definition}
Let $a\in (A,\varphi )$ be a $B$-valued random variable. We say that this
random variable $a$ is $B$-even if

\strut 

\begin{center}
$\varphi (ab_{2}a...b_{m}a)=0_{B},$ whenever $m$ is odd,
\end{center}

\strut 

where $b_{2},...,b_{m}\in B$ are arbitrary. In particular, if $a$ is $B$%
-even, then $\varphi (a^{m})=0_{B},$ whenever $m$ is odd. But the converse
is not true, in general.
\end{definition}

\strut

Recall that in the $*$-probability space model, the $B$-evenness guarantees
the self-adjointness (See [8]). But the above definition is more general. By
using the M\"{o}bius inversion, we have the following characterization ;

\strut

\begin{proposition}
Let $a\in (A,\varphi )$ be a $B$-valued random variable. Then $a$ is $B$%
-even if and only if

\strut 

\begin{center}
$k_{m}\left( \underset{m-times}{\underbrace{a,.......,a}}\right) =0_{B},$
whenever $m$ is odd.
\end{center}
\end{proposition}

\strut

\begin{proof}
($\Rightarrow $) Suppose that $a\in (A,\varphi )$ is $B$-even. Assume that $%
n\in \Bbb{N}$ is odd. Then

\strut

$\ \ \ k_{n}\left( \underset{n-times}{\underbrace{a,.......,a}}\right)
=c^{(n)}\left( a\otimes b_{2}a\otimes ...\otimes b_{n}a\right) $

\strut

where $b_{2},...,b_{n}\in B$ are arbitrary

\strut

$\ \ \ \ \ \ \ =\underset{\pi \in NC(n)}{\sum }\widehat{\varphi }(\pi
)(a\otimes b_{2}a\otimes ...\otimes b_{n}a)\mu (\pi ,1_{n})$

\strut

$\ \ \ \ \ \ \ =0_{B},$

\strut

since every partition $\pi $ contains at least one odd block.

\strut

($\Leftarrow $) Assume that every odd $B$-valued cumulnats of $a\in
(A,\varphi )$ vanishs also assume that $n\in \Bbb{N}$ is odd. Then

\strut

$\varphi \left( ab_{2}a...b_{n}a\right) =\underset{\pi \in NC(n)}{\sum }%
\widehat{c}(\pi )(a\otimes b_{2}a\otimes ...\otimes b_{n}a)=0_{B},$

\strut

since each $\pi $ contains an odd block.
\end{proof}

\strut

The above proposition says that $B$-evenness is easy to veryfy when we are
dealing with either $B$-moments or $B$-cumulants. Now, define a subset $%
NC^{(even)}(2k)$ of $NC(2k),$ for any $k\in \Bbb{N}$ ;

\strut

\begin{center}
$NC^{(even)}(2k)=\{\pi \in NC(2k):\pi $ does not contain odd blocks$\}.$
\end{center}

\strut

We have that ;

\strut

\begin{proposition}
Let $k\in \Bbb{N}$ and let $a\in (A,\varphi )$ be $B$-even. Then

\strut 

\begin{center}
$k_{2k}\left( \underset{2k-times}{\underbrace{a,.......,a}}\right) =%
\underset{\pi \in NC^{(even)}(2k)}{\sum }\widehat{\varphi }(\pi )\left(
a\otimes b_{2}a\otimes ...\otimes b_{2k}a\right) \mu (\pi ,1_{2k})$
\end{center}

\strut 

equivalently,

\strut 

\begin{center}
$\varphi \left( ab_{2}a...b_{2k}a\right) =\underset{\pi \in NC^{(even)}(2k)}{%
\sum }\widehat{c}(\pi )\left( a\otimes b_{2}a\otimes ...\otimes
b_{2k}a\right) .$
\end{center}
\end{proposition}

\strut

\begin{proof}
By the previous proposition, it is enough to show one of the above two
formuli. Fix $k\in \Bbb{N}.$ Then

\strut

$\ \ k_{2k}\left( a,...,a\right) =c^{(2k)}\left( a\otimes b_{2}a\otimes
...\otimes b_{2k}a\right) $

\strut

$\ \ \ \ \ \ \ \ \ \ \ \ \ \ \ \ \ \ \ \ =\underset{\pi \in NC(2k)}{\sum }%
\widehat{\varphi }(\pi )\left( a\otimes b_{2}a\otimes ...\otimes
b_{2k}a\right) \mu (\pi .1_{2k}).$

\strut

Now, suppose that $\theta \in NC(2k)$ and $\theta $ contains its odd block $%
V_{o}\in \pi (o)\cup \pi (i).$ Then

\strut

\strut (2.2.1)

\begin{center}
$\widehat{\varphi }(\theta )\left( a\otimes b_{2}a\otimes ...\otimes
b_{2k}a\right) =0_{B}.$
\end{center}

\strut

Define

\strut

\begin{center}
$NC^{(odd)}(2k)=\{\pi \in NC(2k):\pi $ contains at least one odd block$\}.$
\end{center}

\strut

Then, for any $\theta \in NC^{(odd)}(2k),$ the formular (2.2.1) holds. So,

\strut

\begin{center}
$k_{2k}(a,...,a)=\underset{\pi \in NC(2k)\,\,\setminus \,\,NC^{(odd)}(2k)}{%
\sum }\widehat{\varphi }(\pi )(a\otimes b_{2}a\otimes ...\otimes b_{2k}a)\mu
(\pi ,1_{2k}).$
\end{center}

\strut

It is easy to see that, by definition,

\strut

\begin{center}
$NC^{(even)}(2k)=NC(2k)\,\setminus \,NC^{(odd)}(2k).$
\end{center}
\end{proof}

\strut

\begin{proposition}
Let $a_{1}$ and $a_{2}$ be $B$-even elements in $(A,\varphi ).$ If $a_{1}$
and $a_{2}$ are free over $B,$ then $a_{1}+a_{2}\in (A,\varphi )$ is $B$%
-even, again.
\end{proposition}

\strut

\begin{proof}
Suppose that $a_{1}$ and $a_{2}$ are $B$-free $B$-even elements in $%
(A,\varphi ).$ Let $n\in \Bbb{N}$ be odd. Then

\strut

$\ \ k_{n}\left( (a_{1}+a_{2}),...,(a_{1}+a_{2})\right) $

\strut

$\ \ \ \ \ \ \ \ \ \ \ =k_{n}\left( \underset{n-times}{\underbrace{%
a_{1},......,a_{1}}}\right) +k_{n}\left( \underset{n-times}{\underbrace{%
a_{2},.....,a_{2}}}\right) $

\strut

by $B$-freeness of $a_{1}$ and $a_{2}$

\strut

$\ \ \ \ \ \ \ \ \ \ \ =0_{B}$

\strut

by $B$-evenness of $a_{1}$ and $a_{2}.$ Therefore, by Proposition 2.1, $%
a_{1}+a_{2}$ is also a $B$-even element.
\end{proof}

\strut

Trivially, if $a\in (A,\varphi )$ is $B$-even, then $ba\in (A,\varphi )$ is $%
B$-even, for all $b\in B,$ since we have that

\strut

\begin{center}
$
\begin{array}{ll}
k_{n}\left( \underset{n-times}{\underbrace{ba,.......,ba}}\right) & 
=c^{(n)}\left( ba\otimes b_{2}ba\otimes ...\otimes b_{n}ba\right) \\ 
& =b\cdot c^{(n)}\left( a\otimes b_{2}^{\prime }a\otimes ...\otimes
b_{n}^{\prime }a\right) ,
\end{array}
$
\end{center}

\strut

for all $n\in \Bbb{N},$ where $b_{2},...,b_{n}\in B$ are arbitrary and

\strut

\begin{center}
$b_{2}^{\prime }=b_{2}b,$ \ .... \ , $b_{n}^{\prime }=b_{n}b.$
\end{center}

\strut

\strut

\subsection{$B$-valued R-diagonal Pairs}

\strut

\strut

In this section, we will discuss about $B$-valued R-diagonality of pairs of $%
B$-valued random variables. Likewise, let $(A,\varphi )$ be a NCPSpace over $%
B.$ Remark that to define R-doagonal pairs, we need to assume that $\varphi $
is a $B$-trace.

\strut

\begin{definition}
Let $(A,\varphi )$ be a NCPSpace over $B$ and let $x_{1},x_{2}\in (A,\varphi
)$ be $B$-valued random variables. We say that a pair $(x_{1},x_{2})$ is
R-diagonal if there exists a $B$-formal series $f_{1},\,g\in \Theta _{B}^{1}$
such that

\strut 

\begin{center}
$R_{x_{1},x_{2}}(z_{1},z_{2})=f(z_{1}z_{2})+g(z_{2}z_{1}).$
\end{center}

\strut 

We call this $B$-formal series $\ (f,$ $g)\in \Theta _{B}^{1}\times \Theta
_{B}^{1}$ \ the determining series of the pair $(x_{1},x_{2}).$
\end{definition}

\strut

\begin{theorem}
Let $(A,\varphi )$ be a NCPSpace over $B$ and let $x,y\in (A,\varphi )$ be $B
$-valued random variables. Suppose that the pair $(x,y)$ is an R-diagonal
pair with its determining series \ $(f,g)\in \Theta _{B}^{1}.$ Then

\strut 

\begin{center}
$f=R_{xy}\,\,\frame{*}_{B}\,\,Mob$ \ \ and \ \ $g=R_{yx}\,\,\frame{*}%
_{B}\,\,\,Mob$
\end{center}
\end{theorem}

\strut

\begin{proof}
Let $(x,y)$ be an R-diagonal pair with its determining series $\ (f,g)\in
\Theta _{B}^{1}.$ Then by definition,

\strut

\begin{center}
$R_{x,y}(z_{1},z_{2})=f(z_{1}z_{2})+g(z_{2}z_{1}).$
\end{center}

\strut

Now put

\begin{center}
\strut $f(z)=\sum_{n=1}^{\infty }b_{n}\,z^{n}$ \ and \ $g=\sum_{n=1}^{\infty
}b_{n}^{\prime }\,z^{n}.$
\end{center}

Then

\strut

(2.4.1)

\begin{center}
$R_{x,y}(z_{1},z_{2})=\sum_{n=1}^{\infty
}b_{n}(z_{1}z_{2})^{n}+\sum_{n=1}^{\infty }b_{n}^{\prime
}(z_{2}z_{1})^{n}\in \Theta _{B}^{2}.$
\end{center}

\strut

While, by definition,

\strut

\strut (2.4.2)

\begin{center}
$R_{x,y}(z_{1},z_{2})=\sum_{n=1}^{\infty }\underset{i_{1},...,i_{n}\in
\{1,2\}^{n}}{\sum }k_{n}\left( x_{i_{1}},...,x_{i_{n}}\right)
z_{i_{1}}...z_{i_{n}},$
\end{center}

\strut

\strut where $x_{i_{1}},...,x_{i_{n}}\in \{x,y\},$ for all $%
(i_{1},...,i_{n})\in \{1,2\}^{n},$ $n\in \Bbb{N}.$ By (2.4.1) and (2.4.2),
we can conclude that the R-diagonality of the pair $(x,y)$ makes that the
only nonvanishing mixed $(i_{1},...,i_{m})$-th cumulants of $x$ and $y$
appear when $m$ is even and

\strut

\strut \strut (2.4.3)

\begin{center}
$(i_{1},...,i_{m})=\left( \underset{\frac{m}{2}-times}{\underbrace{%
(1,2),...,(1,2)}}\right) =(1,2,1,2,...,1,2)$
\end{center}

or

\begin{center}
$(i_{1},...,i_{m})=\left( \underset{\frac{m}{2}-times}{\underbrace{%
(2,1),...,(2,1)}}\right) =(1,2,1,2,...,1,2)$
\end{center}

\strut

Therefore, by (2.4.3), we have that the formular (2.4.2) goes to

\strut

(2.4.4)

$\ \ \ \ \ \ \ \ \ =\sum_{n=1}^{\infty }k_{2n}\left( x,y,...,x,y\right)
(z_{1}z_{2})^{n}+\sum_{n=1}^{\infty }k_{2n}\left( y,x,...,y,x\right)
(z_{2}z_{1})^{n}.$

\strut

i.e

\begin{center}
$f(z)=\sum_{n=1}^{\infty }k_{2n}\left( x,y,...,x,y\right) z^{n}$
\end{center}

and

\begin{center}
$g(z)=\sum_{n=1}^{\infty }k_{2n}\left( y,x,...,y,x\right) z^{n},$
\end{center}

in $\Theta _{B}^{1}$

\strut

Observe that

\strut

\begin{center}
$coef_{2n}(f)=coef_{n}\left( R_{xy}\,\,\frame{*}_{B}\,\,Mob\right) $
\end{center}

amd

\begin{center}
$coef_{2n}(g)=coef_{n}\left( R_{yx}\,\,\,\frame{*}_{B}\,\,\,Mob\right) ,$
\end{center}

\strut

for all $n\in \Bbb{N}.$ We will only consider the first case. Fix $n\in \Bbb{%
N}.$ Then

\strut

$\ coef_{n}\left( R_{xy}\right) =k_{n}\left( \underset{n-times}{\underbrace{%
xy,.......,xy}}\right) $

\strut

$\ \ \ \ =c^{(n)}\left( xy\otimes b_{2}xy\otimes ...\otimes b_{n}xy\right) $

\strut

where $b_{2},...,b_{n}\in B$ are arbitrary

\strut

$\ \ \ \ =\underset{\pi \in NC(n)}{\sum }\widehat{\varphi }(\pi )\left(
xy\otimes b_{2}xy\otimes ...\otimes b_{n}xy\right) \mu (\pi ,1_{n})$

\strut

$\ \ \ \ =\underset{\theta \in NC(2n),\,\theta \vee \theta _{0}=1_{2n}}{\sum 
}\widehat{\varphi }(\theta )\left( x\otimes y\otimes b_{2}x\otimes y\otimes
...\otimes b_{n}x\otimes y\right) \mu (\pi ,1_{2n})$

\strut

where $\theta _{0}=\{(1,2),(3,4),...,(2n-1,2n)\}\in NC(2n)$

\strut

\strut (2.4.5)

$\ \ \ \ =k_{2n}\left( x,y,...,x,y\right) ,$

\strut

by (2.4.3). Notice that, in general, the last equality (2.4.5) of the above
formular does not hold true. But, since we have the relation (2,4,3), under
the R-diagonality of $(x,y),$ it holds true. Since

\strut

\begin{center}
$k_{2n}\left( x,y,...,x,y\right) =k_{2n}\left( x,y,...,x,y\right) \cdot
1_{B},$
\end{center}

\strut

(2.4.5) is same as

\begin{center}
\strut $coef_{n}\left( f\,\,\,\frame{*}_{B}\,\,\,Zeta\right) ,$
\end{center}

by the Section 1.2.

\strut

Similarly, we have that, for any fixed $n\in \Bbb{N},$

\strut

\begin{center}
$coef_{n}\left( R_{yx}(z)\right) =coef_{n}\left( g\,\,\,\frame{*}%
_{B}\,\,\,Zeta\right) .$
\end{center}

\strut

Thus

\strut

$R_{xy}=f\,\,\,\frame{*}_{B}\,\,\,Zeta$ \ \ and \ \ $R_{yx}=g\,\,\,\frame{*}%
_{B}\,\,\,Zeta$

\strut

and hence, equivalently,

\strut

\begin{center}
$f=R_{xy}\,\,\,\frame{*}_{B}\,\,\,Mob$ \ \ and \ \ $g=R_{yx}\,\,\,\frame{*}%
_{B}\,\,\,Mob.$
\end{center}
\end{proof}

\strut \strut

We have the following characterization of R-diagonal pairs with respect to $%
B $-valued cumulants ;

\strut

\begin{theorem}
Let $(A,\varphi )$ be a NCPSpace over $B$ and let $x,y\in (A,\varphi )$ be $B
$-valued random variables. Then the pair $(x,y)$ is an R-diagonal pair if
and only if the only nonvanishing mixed $B$-cumulants of $x$ and $y$ are

\strut \strut 

(2.5.1) $\ \ \ \ \ \ \ \ \ \ \ \ \ \ \ \ \ k_{2n}$\strut $\left( \underset{%
2n-times}{\underbrace{x,y,x,y,...,x,y}}\right) $

and

(2.5.2) $\ \ \ \ \ \ \ \ \ \ \ \ \ \ \ \ \ k_{2n}\left( \underset{2n-times}{%
\underbrace{y,x,y,x,...,y,x}}\right) ,$

\strut 

for all $n\in \Bbb{N}.$
\end{theorem}

\strut

\begin{proof}
($\Rightarrow $) By the previous theorem, if the pair of $B$-valued random
variables $(x,y)$ is $R$-diagonal, then

\strut

\begin{center}
$R_{x,y}(z_{1},z_{2})=f(z_{1}z_{2})+g(z_{2}z_{1})\in \Theta _{B}^{2},$
\end{center}

\strut

where $(f,g)\in \Theta _{B}^{1}\times \Theta _{B}^{1}$ is the determining
series of $(x,y)$ such that

\strut

\begin{center}
$f(z)=\left( R_{xy}\,\,\,\frame{*}_{B}\,\,\,Moz\right) (z)$ \ and \ $%
g(z)=\left( R_{yx}\,\,\frame{*}_{B}\,\,Mob\right) (z).$
\end{center}

\strut

By the relation (2.4.3) in the proof of the previous theorem, we can get
that if $(x,y)$ is R-diagonal, then the only nonvanishing mixed cumulants of 
$x$ and $y$ have the form of (2.5.1) or (2.5.2), with respect to the
coefficients of $\ f$ \ and \ $g$ \ appeared in (2.4.4).

\strut

($\Leftarrow $) Conversely, assume that the $B$-valued random variables $x$
and $y$ have their mixed cumulants satisfying that the only nonvanising
mixed cumulants have the form of (2.5.1) or (2.5.2). Then we can easily
construct $f$ and $g$ like (2.4.4). Then the pair $(x,y)$ satisfies that

\strut

\begin{center}
$R_{x,y}(z_{1},z_{2})=f(z_{1}z_{2})+g(z_{2}z_{1}).$
\end{center}

\strut

Therefore, $(x,y)$ is R-diagonal.
\end{proof}

\strut

The following theorem plays a key role for observing the R-transforms of
commutators ;

\strut

\begin{lemma}
Let $a,a^{\prime }\in (A,\varphi )$ be $B$-even. If $a$ and $a^{\prime }$
are free over $B,$ then $\varphi (aa^{\prime })=0_{B}=\varphi (a^{\prime }a).
$
\end{lemma}

\strut

\begin{proof}
Clearly, by the M\"{o}bius inversion, we have that

\strut

$\ \ \ \varphi (aa^{\prime })$\strut $=k_{2}(a,a^{\prime
})+k_{1}(a)k_{1}(a^{\prime })$

$\ \ \ \ \ \ \ \ \ \ \ \ \ =0_{B}+\varphi (a)\varphi (a^{\prime })$

\strut

by the $B$-freeness of $a$ and $a^{\prime }$

\strut

$\ \ \ \ \ \ \ \ \ \ \ \ \ =0_{B}+0_{B}\cdot 0_{B}=0_{B}$

\strut

by the $B$-evenness of $a$ and $a^{\prime }.$ Similarly,

\strut

\begin{center}
$\varphi (a^{\prime }a)=0_{B}.$
\end{center}
\end{proof}

\strut \strut

\begin{theorem}
Let $a,a^{\prime }\in (A,\varphi )$ be $B$-even. If they are free over $B,$
then the pair $(aa^{\prime },a^{\prime }a)$ is $B$-valued R-diagonal.
\end{theorem}

\strut

\begin{proof}
Suppose that $B$-valued random variables $a$ and $a^{\prime }$ are $B$-free $%
B$-even random variables. It suffices to show that the only nonvanishing
mixed cumulants of $aa^{\prime }$ and $a^{\prime }a$ have the form

\strut

\begin{center}
$k_{2n}\left( aa^{\prime },a^{\prime }a,aa^{\prime },a^{\prime
}a,...,aa^{\prime },a^{\prime }a\right) $
\end{center}

or

\begin{center}
$k_{2n}\left( a^{\prime }a,aa^{\prime },a^{\prime }a,aa^{\prime
},...,a^{\prime }a,aa^{\prime }\right) ,$
\end{center}

\strut

for all $n\in \Bbb{N}.$ Put $x=aa^{\prime }$ and $y=a^{\prime }a.$

\strut

Now, fix $n\in \Bbb{N}.$ Suppose that the mixed index $(i_{1},...,i_{n})\in
\{1,2\}^{n}$ is not alternating (i.e, neither $(i_{1},...,i_{n})\neq
(1,2,1,2,...i_{n})$ nor $(i_{1},...,i_{n})\neq (2,1,2,1,...,i_{n})$). Then
we may assume that there exists at least one $j\in \{1,...,n-1\}$ such that $%
i_{j}$ satisfies either $i_{j}=1=i_{j+1}$ or $i_{j}=2=i_{j+1}.$ Let's assume
that $i_{j}=1=i_{j+1}.$ Then

\strut

$k_{n}\left( x_{i_{1}},...x_{i_{j-1}},\underset{j-th}{x},\,\,\underset{j+1-th%
}{x},\,x_{i_{j+2}},...,x_{i_{n}}\right) $

\strut

where $x_{i_{1}},...,x_{i_{n}}\in \{x,y\}$

\strut

\strut (2.7.1)

\strut

$\ \ =\underset{\pi \in NC(n)}{\sum }\widehat{\varphi }(\pi
)(x_{i_{1}}\otimes b_{i_{2}}x_{i_{2}}\otimes ...\otimes
b_{i_{j-1}}x_{i_{j-1}}\otimes b_{i_{j}}x\otimes b_{i_{j+1}}x$

\begin{center}
$\otimes b_{i_{j+2}}x_{i_{j+2}}\otimes ...\otimes b_{i_{n}}x_{i_{n}})\mu
(\pi ,1_{n}),$

\strut
\end{center}

where $b_{i_{2}},...,b_{i_{n}}\in B$ are arbitrary.

\strut

First, observe that $\varphi (x)=0_{B}=\varphi (y)$ ;

\strut

\begin{center}
$\varphi (x)=\varphi (aa^{\prime })=0_{B}=\varphi (a^{\prime }a)=\varphi
(y), $
\end{center}

\strut

by the previous lemma. Therefore, for any partitions, $\theta ,$ in $NC(n)$
containing singleton blocks $(j)$ and $(j+1)$, $\widehat{\varphi }(\theta
)(...)$ vanish. So, the formular (2.7.1) is same as

\strut

$\ \ \underset{\pi \in S}{\sum }\widehat{\varphi }(\pi )(x_{i_{1}}\otimes
b_{i_{2}}x_{i_{2}}\otimes ...\otimes b_{i_{j-1}}x_{i_{j-1}}\otimes
b_{i_{j}}x\otimes b_{i_{j+1}}x$

\begin{center}
$\otimes b_{i_{j+2}}x_{i_{j+2}}\otimes ...\otimes b_{i_{n}}x_{i_{n}})\mu
(\pi ,1_{n}),$
\end{center}

\strut

where

\strut

\begin{center}
$S=\{\pi \in NC(n):(j)\notin \pi \,\,\,\,\&\,\,\,(j+1)\notin \pi \}.$
\end{center}

\strut

Second, observe that $\varphi (xbx)=0_{B},$ for any $b\in B$ ;

\strut

$\varphi (xbx)=\varphi (aa^{\prime }baa^{\prime })=\underset{\pi \in NC(4)}{%
\sum }\widehat{c}(\pi )\left( a\otimes a^{\prime }\otimes ba\otimes
a^{\prime }\right) $

\strut

$\ \ \ \ \ \ \ \ \ \ =\widehat{c}(0_{4})\left( a\otimes a^{\prime }\otimes
ba\otimes a^{\prime }\right) $

\strut

by the $B$-freeness of $a$ and $a^{\prime }$ and by the $B$-evenness of $a$
and $a^{\prime }$

\strut

$\ \ \ \ \ \ \ \ \ \ =k_{1}(a)\cdot k_{1}(a^{\prime })\cdot \left(
bk_{1}(a)\right) \cdot k_{1}(a^{\prime })$

\strut

$\ \ \ \ \ \ \ \ \ \ =\varphi (a)\cdot \varphi (a^{\prime })\cdot b\varphi
(a)\cdot \varphi (a^{\prime })=0_{B},$

\strut

by the $B$-evenness of $a$ and $a^{\prime }.$ So, this shows that the
formular (2.7.1) goes to

\strut

(2.7.2)

$\ \ \ \ \ \ \ \ \underset{\pi \in S^{\prime }}{\sum }\widehat{\varphi }(\pi
)(x_{i_{1}}\otimes b_{i_{2}}x_{i_{2}}\otimes ...\otimes
b_{i_{j-1}}x_{i_{j-1}}\otimes b_{i_{j}}x\otimes b_{i_{j+1}}x$

\begin{center}
$\otimes b_{i_{j+2}}x_{i_{j+2}}\otimes ...\otimes b_{i_{n}}x_{i_{n}})\mu
(\pi ,1_{n}),$
\end{center}

\strut

where

\strut

\begin{center}
$S^{\prime }=\{\pi \in NC(n):\pi $ does not contain $(j),(j+1),(j,j+1)\}.$
\end{center}

\strut

Consider the set $S^{\prime }.$ Suppose that there exists at least one $%
k\neq j$ in $\{1,...,n-1\}$ such that $i_{k}$ satisfies either $%
x_{i_{k}}=x=x_{i_{k+1}}$ or $x_{i_{k}}=y=x_{i_{k+1}}.$ Then we can do the
same job as before on $S^{\prime }$ and we can get a set $S^{\prime \prime
}. $ Inductively we have that $S^{(p)}$ such that

\strut

$\ k_{n}\left( x_{i_{1}},...,x_{i_{n}}\right) $

\strut

\strut (2.7.3)

$\ \ \ =\underset{\pi \in S^{(p)}}{\sum }\widehat{\varphi }(\pi
)(x_{i_{1}}\otimes b_{i_{2}}x_{i_{2}}\otimes ...\otimes
b_{i_{j-1}}x_{i_{j-1}}\otimes b_{i_{j}}x\otimes b_{i_{j+1}}x$

\begin{center}
$\otimes b_{i_{j+2}}x_{i_{j+2}}\otimes ...\otimes b_{i_{n}}x_{i_{n}})\mu
(\pi ,1_{n}),$
\end{center}

\strut \strut

where $S^{(p)}$ is determined by the $p$-induction of the previous process.
Now, let's assume that the formular (2.7.3) does not vanish. Then the mixed
index $(i_{1},...,i_{n})\in \{1,2\}^{n}$ should be alternating. But it
contradict our assumption.

\strut

We can get the same result, when we replace $x_{i_{j}}=x=x_{j+1}$ by $%
x_{i_{j}}=y=x_{i_{j+1}}.$

\strut

Now, we have to observe that $n$ should be even. Suppose that $n$ is odd and
we have an alternating mixed index $(1,2,1,2,...,1,2,1).$ Then

\strut

$k_{n}\left( x,y,x,y,...,x,y,x\right) $

\strut

$\ \ =\underset{\pi \in NC(n)}{\sum }\widehat{\varphi }(\pi )(x\otimes
b_{2}y\otimes b_{3}x\otimes b_{4}y\otimes $

\begin{center}
$...b_{n-2}x\otimes b_{n-1}y\otimes b_{n}x)\mu (\pi ,1_{n}),$
\end{center}

\strut

where $b_{2},...,b_{n}\in B$ are arbitrary

\strut

$\ \ =\underset{\pi \in NC(2n),\,\pi \vee \theta =1_{2n}}{\sum }\widehat{%
\varphi }(\pi )(a\otimes a^{\prime }\otimes b_{2}a^{\prime }\otimes a\otimes 
$

\begin{center}
$...\otimes b_{n-2}a\otimes a^{\prime }\otimes b_{n-1}a^{\prime }\otimes
a\otimes b_{n}a\otimes a^{\prime })\mu (\pi ,1_{n}),$
\end{center}

\strut

where

\begin{center}
$\theta =\{(1,2),...,(1,2)\}\in NC(2n),$

\strut

$\ \ =0_{B},$

\strut

by [15]. Similarly, we can get the same result if we replace $(1,2,...,1,2)$
by $(2,1,...,2,1).$ Therefore, by Theorem 2.5, the pair $(x,y)=(aa^{\prime
},a^{\prime }a)$ is an $R$-diagonal pair. Similarly, we can conclude that
the pair $(y,x)=(a^{\prime }a,aa^{\prime })$ is R-diagonal.
\end{center}
\end{proof}

\strut

\strut \strut

\subsection{\strut $B$-valued R-diagonal Elements in $C^{*}$-Probability
Spaces over $B$}

\strut

\strut

\strut In this chapter, we will consider the $B$-valued R-diagonal elements
in a $C^{*}$-probability space over a unital $C^{*}$-algebra $B,$ where the $%
B$-functional $\varphi :A\rightarrow B$ is a normalized positive $B$%
-functional. We say that a $B$-valued random variable in $(A,\varphi )$ is $%
B $-even if it is a self-adjoint $B$-even element (in the sense of Section
2.1). Let $x\in (A,\varphi )$ be a $B$-valued random variable. We say that $%
x\in (A,\varphi )$ is R-diagonal if the pair $\left( x,x^{*}\right) $ is an
R-diagonal pair. By the characterization considered in Section 2.2, we can
redefine that a $B$-valued random variable $x$ is R-diagonal if the only
nonvanishing mixed cumulants of $x$ and $x^{*}$ are of the form

\strut

\begin{center}
$k_{2n}\left( x,x^{*},...,x,x^{*}\right) $ \ and \ $k_{2n}\left(
x^{*},x,...,x^{*},x\right) ,$
\end{center}

\strut

for all $n\in \Bbb{N}.$\strut

\strut

\strut

\strut \textbf{References}

\strut 

\strut 

{\small [1] \  A. Nica, R-transform in Free Probability, IHP course note,
available at www.math.uwaterloo.ca/\symbol{126}anica.}

{\small [2] \ \ A. Nica, R-transforms of Free Joint Distributions and
Non-crossing Partitions, J. of Func. Anal, 135 (1996), 271-296.}

{\small [3] \ \ A. Nica and R. Speicher, R-diagonal Pair-A Common Approach
to Haar Unitaries and Circular Elements, (1995), www.mast.queensu.ca/\symbol{%
126}speicher.}

{\small [4] \ \ A. Nica, R-diagonal Pairs Arising as Free Off-diagonal
Compressions, available at www.math.uwaterloo.ca/\symbol{126}anica.}

{\small [5] \ \ A. Nica, D. Shlyakhtenko and R. Speicher, R-diagonal
Elements and Freeness with Amalgamation, Canad. J. Math. Vol 53, Num 2,
(2001) 355-381.}

{\small [6] \ \ A. Nica and R.Speicher, Commutators of Free Random
Variables, Duke Math J, Vol. 92, No. 3 (1998) 553 - 392.}

{\small [7] \ \ D.Voiculescu, K. Dykemma and A. Nica, Free Random Variables,
CRM Monograph Series Vol 1 (1992).}

{\small [8] \ \ I. Cho, Amalgamated Boxed Convolution and Amalgamated
R-transform Theory (2002), Preprint.}

{\small [9] \ \ I. Cho, I. Cho, R-transform Theory of Commutators of
Amalgamated Random Variables (2004) (preprint).}

{\small [10] R. Speicher Combinatorial Theory of the Free Product with
Amalgamation and Operator-Valued Free Probability Theory, AMS Mem, Vol 132 ,
Num 627 , 1998.}

{\small [11] R. Speicher, Combinatorics of Free Probability Theory IHP
course note, available at www.mast.queensu.ca/\symbol{126}speicher.}

\label{REF}

\strut

\end{document}